\documentclass[twoside]{article}

\usepackage{psfrag}
\usepackage{graphicx}
\usepackage{amssymb}

\newtheorem{theorem}{Theorem}
\newtheorem{proposition}{Proposition}
\newtheorem{definition}{Definition}
\newtheorem{lemma}{Lemma} 
\newtheorem{corollary}{Corollary}

\def\zf{I\! I}

\begin{document}

\date{}

\title{Hyperbolic Rank of Products}
\maketitle
\begin{center}
{\large Thomas Foertsch* $\hspace{1cm}$ and $\hspace{1cm}$ Viktor Schroeder}
\footnote{* supported by SNF Grant 21 - 589 38.99}
\end{center}


\vspace{0.5cm}

\begin{abstract}
Generalizing \cite{brfa} we prove
the existence of a bilipschitz embedded manifold of pinched negative 
curvature and dimension $m_1+m_2-1$ in the product $X:=X_1^{m_1}\times X_2^{m_2}$ of
two Hadamard manifolds $X_i^{m_i}$ of dimension $m_i$ with pinched negative curvature. \\
Combining this result with \cite{buysch} we prove the additivity of the hyperbolic rank for products of manifolds with pinched 
negative curvature.
\end{abstract}


\section{Introduction}

In \cite{brfa} the authors proved that the product $X:=H^{m_1}\times ... \times H^{m_k}$ of hyperbolic spaces $H^{m_i}$
admits an embedding of $Y:=H^{m_1+...+m_k-k+1}$ in $X$ that is quasi-isometric in the sense that the Riemannian 
distance functions $d_X$ and $d_Y$ on $X$ and $Y$ are related via
\begin{displaymath}
d_X|_Y \; \le d_Y \le \; \alpha \cdot d_X|_Y \; + \; \beta \;\;\; ,
\end{displaymath}
with constants $\alpha , \beta \in \mathbb{R}^+$, $\alpha > 1$. \\
In \cite{foe} this was generalized to the case of Riemannian products of
certain types of warped products, that were shown to admit bilipschitz 
embeddings of warped products of negative sectional curvature. \\
In this paper we further generalize this observation proving the 
\begin{theorem} \label{prodhadamard}
Let $(X_i^{m_i},g_i)$, $i=1,...,k$, be Hadamard manifolds of pinched negative 
sectional curvature $-b_i^2 \le K(X_i^{m_i},g_i) \le -a_i^2$. Then their 
Riemannian product $(X={\Pi}_{i=1}^k X_i^{m_i} ,g_X)$ admits a bilipschitz embedding
of a Hadamard manifold $(Y^n,g_Y)$ of pinched negative sectional curvature 
$-b^2\le K(Y^n,g_Y)\le -a^2$ and dimension $n:=m_1+...+m_k-k+1$.
\end{theorem}

The existence result of Theorem \ref{prodhadamard} implies the additivity of the hyperbolic rank for products
of Hadamard manifolds with pinched negative curvature (compare \cite{buysch}). \\
Given a metric space $X$ consider all locally compact $CAT(-1)$ Hadamard spaces $Y$ quasi-isometrically embedded
into $X$ and let
\begin{displaymath}
rank_h X \; := \; \sup\limits_{Y} \; dim \, {\partial}_{\infty} Y
\end{displaymath}
over all such $Y$, where $dim {\partial}_{\infty}Y$ is the topological dimension of the boundary of $Y$.
This invariant was introduced by Gromov (\cite{gr}) in a slightly stronger form, where he called it the hyperbolic corank.
As explained in \cite{buysch} we will call it the hyperbolic rank. \\
Combining Theorem \ref{prodhadamard} with the result of \cite{buysch} we obtain
\begin{theorem} \label{theo-hyprank}
Let $X=X_1\times ... \times X_k$ be the product of Hadamard manifolds with pinched negative curvature, then 
\begin{displaymath}
rank_h X \; = \; \sum\limits_{i=1}^k \, rank_h X_i .
\end{displaymath}
\end{theorem}


\section{Proof of Theorem \ref{prodhadamard}}

It is convinient to have the following 
\begin{definition}
Let $(X,g_X)$ be a Riemannian manifold. A function $f:X\longrightarrow \mathbb{R}$ is said to have pinched convex levels,
if there are constants $0<c_1\le c_2$ such that for any unit speed geodesic $c:(-\eta ,\eta)\longrightarrow Y$ with
$\dot{c}(0)\perp grad f(c(0))$:
\begin{displaymath}
c_1 \; \le \; (f\circ c)''(0) \; \le \; c_2 .
\end{displaymath}
\end{definition}

We now prove the 
\begin{lemma}  \label{lemma}
Let $(Y,g_Y)$ be a complete Riemannian manifold with
\begin{description}
\item[1)] $|K|\le \Lambda$, 
\end{description}
such that there exists a function $f:Y\longrightarrow \mathbb{R}$ satisfying
\begin{description}
\item[2)] $|| \, grad \; f \, || \equiv 1$,
\item[3)] $f$ has pinched convex levels and
\item[4)] $K({\sigma}_p )< -\delta$ for all planes ${\sigma}_p$ with $grad \; f(p) \subset {\sigma}_p$.
\end{description}
Then $(Y,g_Y)$ is lipshitz equivalent to a Riemannian manifold of pinched negative sectional curvature.
\end{lemma}

{\bf Proof of Lemma \ref{lemma}:}\\
First note that condition $2)$ implies that $Y$ splits into a product
$Y=\mathbb{R}\times N$ such that
\begin{description}
\item[a)] $\mathbb{R}\times p$ is orthogonal to $N^t := t\times N$ for all $t\in \mathbb{R}$ and $p\in N$ and $\mathbb{R}$ 
lies totally geodesic in $Y$, i.e. locally the metric is given by
\begin{displaymath}
g_{Y} \; = \; 
\left(
\begin{array}{cc}
r & 0  \;\; . \;\; . \;\; . \;\; 0 \\
\begin{array}{c}
0 \\
\stackrel{\cdot}{\stackrel{\cdot}{\cdot}} \\
0
\end{array}
& g_{ij}(t,y^k)
\end{array}
\right), \hspace{1.5cm} r\in\mathbb{R}
\end{displaymath}
\end{description}
Now $1)$ and $3)$ guarantee that
\begin{description}
\item[b)] the Riemannian embedding $N\longrightarrow (Y,g_Y)$ is strictly convex in the sense, that its second fundamental form 
tensor is strictly positive definite, i.e. $\forall$ $B \in T_pN$ with $||B||=1$:   
\begin{displaymath}
\zf (B,B) \;\;\;\; > \; \epsilon \;\;\;\;\; \mbox{and}
\end{displaymath}
\item[c)] the intrinsic sectional curvature's absolute value $|K^N(\Delta)|$ of any 
tangent plane $\Delta$ that is tangent to $N$ is bounded by some 
real number $C_1$: 
\begin{displaymath}
|K^{N}(\Delta)| \; < \; C_1,
\end{displaymath}
\end{description}
From $1)$ and $4)$ it further follows that 
\begin{description}
\item[d)] the sectional curvature  $K(\Sigma)$ of any tangent plane $\Sigma$ spanned by a vector tangent to $N$ and 
the coordinate vector $\frac{\partial}{\partial t}|_t$ is bounded by negative constants 
$- {\delta}',-\delta $:  
\begin{displaymath}
-{\delta}' \; < \; K(\Sigma ) \; < \; -\delta .
\end{displaymath}
\end{description}
Finally $1)$ can be rewritten as
\begin{description}
\item[e)] the norm of the curvature tensor $R^{g_Y}$ of $g_Y$ being bounded:
\begin{displaymath}
\parallel R \parallel \; := \; \sup\limits_{\parallel A_i \parallel \le 1} \Big\{ <R_{A_1\; A_2} A_3, A_4>\Big\} \; < \; C_2.
\end{displaymath}
\end{description}

We now stretch the metric into the direction of $grad \, f$ by some factor. Thus consider locally for $\lambda \neq 0$ the metric
\begin{displaymath}
h^{\lambda} \; = \; 
\left(
\begin{array}{cc}
\frac{r}{{\lambda}^2} & 0  \;\; . \;\; . \;\; . \;\; 0 \\
\begin{array}{c}
0 \\
\stackrel{\cdot}{\stackrel{\cdot}{\cdot}} \\
0
\end{array}
& g_{ij}(t,y^k)
\end{array}
\right).
\end{displaymath}
In the following we indicate all objects coming from the metric $h^{\lambda}$ with an additional subscript $\lambda$. \\
Let $\Gamma$ be an arbitrary tangent plane in $TY$  spanned by the two $h^{\lambda}$-orthonormal vectors
\begin{displaymath}
A+B\; = \; a\frac{\partial}{\partial t}{\Big|}_{(t,y^k)} \; + \; b^i\frac{\partial}{\partial y^i}{\Big|}_{(t,y^k)}    \hspace{1cm}
\mbox{and} \hspace{1cm} C \; = \; c^i\frac{\partial}{\partial y^i}{\Big|}_{(t,y^k)}.
\end{displaymath}
The sectional curvature $K^{\lambda}(A+B,C)$ is given via
\begin{eqnarray}
K^{\lambda}(A+B,C) & = & \;\; <R^{\lambda}_{A \; C}A,C>_{\lambda} \; + \;
2\cdot <R^{\lambda}_{A \; C}B,C>_{\lambda} \nonumber \\ 
& & + \; <R^{\lambda}_{B \; C}B,C>_{\lambda}. \label{schnittkr1}
\end{eqnarray}

Let ${\Gamma}^k_{ij}$ denote the Christoffel symbols of $g_Y$ and $\tilde{\Gamma}^k_{ij}$ those of $h^{\lambda}$.
It is obvious that ${\Gamma}^k_{ij}=\tilde{\Gamma}^k_{ij}$ for $k\neq 1$ and $\tilde{\Gamma}^1_{ij}={\lambda}^2{\Gamma}^1_{ij}$.
Further note that ${\Gamma}^1_{1j}=0=\tilde{\Gamma}^1_{1j}$. 
Write $R_{{\partial}^i{\partial}^j} {\partial}^k=\sum_{l}R_{ijk}^l {\partial}^l$,
thus it follows immediately from
\begin{displaymath}
R^l_{ijk} \; = \; \frac{\partial {\Gamma}^l_{ik}}{\partial y^j} \; - \; \frac{\partial {\Gamma}^l_{jk}}{\partial y^i} \; + \;
{\Gamma}^h_{ik} {\Gamma}^l_{jh} \; - \; {\Gamma}^h_{jk} {\Gamma}^l_{ih} 
\end{displaymath}
that $<R^{\lambda}_{A \; C}A,C>_{\lambda}=<R_{A \; C}A,C>$ and
$<R^{\lambda}_{A \; C}B,C>_{\lambda}=<R_{A \; C}B,C>$. With that, equation (\ref{schnittkr1}) turns into
\begin{eqnarray}
K^{\lambda}(A+B,C) & = & \;\; K(A,C) \; {\lambda}^2 \; <A,A>_{\lambda} \; <C,C>_{\lambda}   \nonumber \\
& & + \; 2 \cdot \lambda \; <R_{\frac{A}{\lambda} \; C}B,C>  \label{schnittkr2} \\ 
& & + \; K^{\lambda}(B,C) \; <B,B>_{\lambda} \; <C,C>_{\lambda}. \nonumber 
\end{eqnarray}
Since $g_Y|_N=h^{\lambda}|_N$ the intrinsic curvature of $N$ is independent of the particular choice of $\lambda$.
The Gauss equation yields with $\tilde{B}=\mu B$, $<\tilde{B},\tilde{B} >_{\lambda}=1$:
\begin{displaymath}
\begin{array}{l}
K^{\lambda}(\tilde{B},C)   \\
= \; K^N(\tilde{B},C) 
\; - \; \Big[ {\zf}^{\lambda}(\tilde{B},\tilde{B}) {\zf}^{\lambda}(C,C)
\; - \; {\zf}^{\lambda}(\tilde{B},C) {\zf}^{\lambda}(\tilde{B},C) \Big]  \\
= \; K^N(\tilde{B},C) 
\; - \; {\lambda}^2 \cdot \Big[ {\zf}(\tilde{B},\tilde{B}) {\zf}(C,C)
\; - \; {\zf}(\tilde{B},C) {\zf}(\tilde{B},C) \Big]  \\
\le \; C_1  \; - \; {\lambda}^2 \cdot {\epsilon}^2 ,
\end{array}
\end{displaymath}
where the inequality follows from $b)$. \\
Thus equation (\ref{schnittkr2}) takes the form
\begin{eqnarray} 
K^{\lambda}(A+B,C) & \le & 
-\delta \cdot {\lambda}^2  \cdot <A,A>_{\lambda} \; <C,C>_{\lambda}   \nonumber \\
& & + \big( C_1 - {\lambda}^2 \cdot {\epsilon}^2 \big) \cdot <B,B>_{\lambda} \; <C,C>_{\lambda}   \nonumber \\
& & + 2 \lambda \cdot C_2 \nonumber \\
& \le & - {\lambda}^2 \min \{ \delta , {\epsilon}^2 \} \; + \; C_1 \; + \; 2 \lambda C_2.
\label{upperbound}
\end{eqnarray}
Hence for $\lambda \ge \frac{C_2+\sqrt{{C_2}^2+C_1 \cdot \min \{ \delta , {\epsilon}^2 \}}}{\min \{ \delta , {\epsilon}^2 \}}$, 
$h^{\lambda}$ makes $(Y,h^{\lambda})$ a Riemannian manifold of negative sectional curvature. \\

$|K^{(Y,g_Y)}|<\Lambda$ and the Gauss equation yield an upper bound ${{\epsilon}'}^2$ for
\begin{displaymath}
\Big| {\zf}(\tilde{B} ,\tilde{B}) {\zf}(C,C)
\; - \; {\zf}(\tilde{B},C) {\zf}(\tilde{B},C) \Big| 
\; \le \; {{\epsilon}'}^2
\end{displaymath}
and analogue to equation (\ref{upperbound}) one finds
\begin{equation}
K^{\lambda}(A+B,C) \; \ge \; 
- {\lambda}^2 \max \{ {\delta}' , {{\epsilon}'}^2 \} \; - \; C_1 \; - \; 2 \lambda C_2.
\end{equation}
Furthermore it is obvious that $h^{\lambda}$ is Lipschitz equivalent to $g_Y$. 
\begin{flushright}
{\bf q.e.d.}
\end{flushright}

Applying Lemma \ref{lemma} we can prove the following

\begin{proposition} \label{prophyprank}
Let  $(X_i,g_{X_i})$, $i=1,2$, be Riemannian manifolds satisfying conditions $1)-4)$ of Lemma \ref{lemma}.
Then the natural Riemannian embedding $(Y,g_Y)$ of $Y=\mathbb{R}\times N_1 \times N_2$ in 
$(X,g_X)=(X_1,g_{X_1})\times (X_2,g_{X_2})$,
with
\begin{displaymath}
\begin{array}{cccc}
i: & (t,y^k,\tilde{y}^k) & \longrightarrow & (t,y^k,t,\tilde{y}^k),
\end{array}
\end{displaymath}    
is Lipschitz equivalent to a Riemannian manifold of pinched negative sectional curvature.
\end{proposition}

{\bf Proof of Proposition \ref{prophyprank}:} \\
As in the proof of Lemma \ref{lemma} we first note that conditions $1)-4)$ imply the following: \\
$X_i$ is a product manifold $X_i=\mathbb{R}\times N_i$ of the real numbers $\mathbb{R}$ and a $n_i$-dimensional 
manifold $N_i$ and the metric $g_{X_i}$ is such that
\begin{description}
\item[i)] $\mathbb{R}\times p$ is orthogonal to $N^t_i := t\times N_i$ for all $t\in \mathbb{R}$ and $p\in N_i$, and $\mathbb{R}$ lies
totally geodesic in $X_i$:
\begin{displaymath}
g_{X_1} \; = \; 
\left(
\begin{array}{cc}
1 & 0  \;\; . \;\; . \;\; . \;\; 0 \\
\begin{array}{c}
0 \\
\stackrel{\cdot}{\stackrel{\cdot}{\cdot}} \\
0
\end{array}
& g_{ij}(t,y^k)
\end{array}
\right), \hspace{1cm}
g_{X_2} \; = \; 
\left(
\begin{array}{cc}
1 & 0  \;\; . \;\; . \;\; . \;\; 0 \\
\begin{array}{c}
0 \\
\stackrel{\cdot}{\stackrel{\cdot}{\cdot}} \\
0
\end{array}
& \tilde{g}_{ij}(t,\tilde{y}^k)
\end{array}
\right).
\end{displaymath}
\item[ii)] the Riemannian embeddings $N_i\longrightarrow (X_i,g_{X_i})$ are strictly convex in the sense, 
that their second fundamental form tensors ${\zf}^i$ are strictly positive definite, i.e. $\forall$ $A_i \in T_pN_i$ with 
$||A_i||_i=1$:   
\begin{displaymath}
{\zf}^i(A_i,A_i) \;\;\;\; > \; {\epsilon}_i \;\; ,
\end{displaymath}
\item[iii)] the intrinsic sectional curvatures' absolute values $|K^{N_i}({\Delta}^i)|$ 
of any tangent plane ${\Delta}^i$ that is tangent to $N^i$ 
is bounded by some real number $C^i_1$: 
\begin{displaymath} 
|K^{N_i}({\Delta}^i)| \; < \; C^i_1,
\end{displaymath}
\item[iv)] the sectional curvature  $K^{x_i}({\Sigma}^i)$ of any tangent plane ${\Sigma}^i$ spanned by a vector tangent to $N_i$ and 
the coordinate vector $\frac{\partial}{\partial t^i}|_{t^i}$ is bounded by negative constants
\begin{displaymath}
-{{\delta}'}_i \; < \; K^{X_i}({\Sigma}^i ) \; <  \; -{\delta}_i,
\end{displaymath}
\item[v)] the norms of the curvature tensors $R^{X_i}$ of $g_{X_i}$ are bounded:
\begin{displaymath}
\parallel R^{X_i} \parallel \; := \; \sup\limits_{\parallel A_j \parallel \le 1} 
\Big\{ R^{X_i}_{A_1\; A_2} A_3 A_4\Big\} \; < \; C^i_2.
\end{displaymath}
\end{description}

We now establish that the Riemannian embedding $(Y,g_Y)$ of $Y=\mathbb{R}\times N_1 \times N_2$ in 
$(X,g_X)=(X_1,g_{X_1})\times (X_2,g_{X_2})$, determined by
\begin{displaymath}
\begin{array}{cccc}
i: & (t,y^k,\tilde{y}^k) & \longrightarrow & (t,y^k,t,\tilde{y}^k),
\end{array}
\end{displaymath}
satisfies the requirements a) - e) of Lemma \ref{lemma}:
\begin{description}
\item[a)] The pullback metric $g_Y:=i^*g_X$ on $Y$ is 
\begin{equation} \label{metricg_Y}
g_{Y} \; = \; 
\left(
\begin{array}{cc}
2 & 0  \;\; . \;\; . \;\; . \;\; 0 \\
\begin{array}{c}
0 \\
\stackrel{\cdot}{\stackrel{\cdot}{\cdot}} \\
0
\end{array}
& g_{N}(t,y^k,\tilde{y}^l)
\end{array}
\right). 
\end{equation}
\item[b)] Let $V=\tilde{B}+\tilde{C}$ be the unique decomposition of an arbitrary vector $V\in T_{(p,q)}N$ into the sum of the
lifts $\tilde{B}$, $\tilde{C}$ of vectors $B\in T_pN_1$ and $C\in T_qN_2$ onto $T_{(p,q)}N$. \\
With ${\zf}^1(B,B)=2{\zf}^Y(\tilde{B},\tilde{B})$, ${\zf}^2(C,C)=2{\zf}^Y(\tilde{C},\tilde{C})$ and \linebreak
${\zf}^Y(\tilde{B},\tilde{C})=0$
it follows
\begin{eqnarray}
\frac{{\zf}^Y(\tilde{B}+\tilde{C},\tilde{B}+\tilde{C})}{g_Y(\tilde{B}+\tilde{C},\tilde{B}+\tilde{C})} 
& = & 
\frac{{\zf}^Y(\tilde{B},\tilde{B}) \; + \; {\zf}^Y(\tilde{C},\tilde{C})}{g_Y(\tilde{B}+\tilde{C},\tilde{B}+\tilde{C})} \nonumber \\
& = & 
\frac{1}{2} \; \frac{{\zf}^1(B,B) \; + \; {\zf}^2(C,C)}{g_Y(\tilde{B}+\tilde{C},\tilde{B}+\tilde{C})} \nonumber \\
& \ge &
\frac{1}{2} \; \min \Big\{ \frac{{\zf}^1(B,B)}{g_{N_1}(B,B)}, \frac{{\zf}^2(C,C)}{g_{N_2}(C,C) } \Big\} \nonumber \\
& = & \frac{1}{2} \; \min \Big\{ {\epsilon}_1,{\epsilon}_2 \Big\} \; =: \; \epsilon . \nonumber
\end{eqnarray}
\item[c)] Consider a tangent plane $\Gamma$ in $TY$ at $\{ t,p,q\} \in Y$ that is tangent to $N^t$. But $N^t$ is the 
Riemannian product $N^t=N_1^t \times N_2^t$. Therefore it is 
\begin{displaymath}
|K^N(\Gamma )| \; \le \; \max \Big\{ C_1^1, \; C_1^2 \Big\} \; =: \; C_1.
\end{displaymath}
\item[d)] Let $\tilde{B},\tilde{C} \; \in \; T_{(t,p,q)}Y$ be as in b), thus tangent to $N$, and such that $\tilde{B}+\tilde{C}$ has unit length
in $(Y,g_Y)$ and set $\tilde{A}:=\frac{1}{\sqrt{2}}\frac{\partial}{\partial t}|_{(t,p,q)}$. Thus $\tilde{A}$ and $\tilde{B}+\tilde{C}$ are orthonormal 
in $(Y,g_Y)$. Recall that $n_i:=dim\, N_i$, $i=1,2$. Coordinates $(t,y^i)$ in $X_1$ 
and $(t,\tilde{y}^k)$ in $X_2$ naturally induce coordinates in $X=X_1\times X_2$ and thus $1+n_1+n_2$ coordinate functions 
$(t,y^i,\tilde{y}^k)$ on $Y$. 
In the following the various indecies are $i,j \in \{ 2,3,...,n_1+1\}$ and $k,l \in \{ n_1+2,....,n_1+n_2+1\}$. With that one 
calculates
\begin{displaymath}
R^Y_{11ij}=2R^{X_1}_{11ij}, \hspace{1cm} R^Y_{11kl}=2R^{X_2}_{11kl} \hspace{1cm} \mbox{and} \hspace{1cm} R^Y_{11ik}=0, 
\end{displaymath}    
and therefore
\begin{eqnarray}
K^Y(\tilde{A},\tilde{B}+\tilde{C}) & = & <R^Y_{\tilde{A}\tilde{B}}\tilde{A},\tilde{B}>_Y \; + \; 
<R^Y_{\tilde{A}\tilde{C}}\tilde{A},\tilde{C}>_Y  \nonumber \\
& = & 2 <R^{X_1}_{A_1B}A_1,B>_{X_1} \; + \;  2 <R^{X_2}_{A_2C}A_2,C>_{X_2} \nonumber \\
& \stackrel{(*)}{=} & K^{X_1}(A_1, B) \, <B,B>_{X_1} \; + \; K^{X_2}(A_2, C) \, <C,C>_{X_2}. \nonumber 
\end{eqnarray}
Note, that in $(*)$ we made use of $<A_i,A_i>_{X_i}=\frac{1}{\sqrt{2}}$, $i=1,2$. \\
With that we finally achieve
\begin{displaymath}
-{\delta}' \; := \; \min \Big\{ -{{\delta}'}_1, \; -{{\delta}'}_2 \Big\} \; \le \;
K^Y(\tilde{A},\tilde{B}+\tilde{C}) \; \le \;
\max \Big\{ -{\delta}_1, \; -{\delta}_2 \Big\} \; =: \; - \delta . \nonumber
\end{displaymath}
\item[e)] In order to show
\begin{displaymath}
\parallel R^Y \parallel \; := \; \sup\limits_{\parallel A_i \parallel \le 1} \Big\{ <R^Y_{A_1\; A_2} A_3, A_4>_Y\Big\} \; < \; C_2
\end{displaymath}
we emphasize that the sectional curvature's absolute value $|K^Y|$ is bounded: \\
Let therefore $\Gamma$ be an arbitrary tangent plane in $TY$ spanned by the two
$g_Y$-orthonormal vectors 
\begin{displaymath}
A+B\; = \; a\frac{\partial}{\partial t}{\Big|}_{(t,y^k)} \; + \; b^i\frac{\partial}{\partial y^i}{\Big|}_{(t,y^k)}    \hspace{1cm}
\mbox{and} \hspace{1cm} C \; = \; c^i\frac{\partial}{\partial y^i}{\Big|}_{(t,y^k)}.
\end{displaymath}
For the sectional curvature of $\Gamma$ one has
\begin{eqnarray}
K^Y(\Gamma ) & = & \;\; \underbrace{<R^Y_{A \; C}A,C>_Y}_{A)} \; + \;
2 \underbrace{<R^Y_{A \; C}B,C>_Y}_{B)} \nonumber \\ 
& & + \; \underbrace{<R^Y_{B \; C}B,C>_Y}_{C)}. \nonumber
\end{eqnarray}
Thus it sufficies to show that the terms $A)$, $B)$ and $C)$ are bounded: \\
\begin{description}
\item[A)] $<R^Y_{A \; C}A,C>_Y$ is bounded due to $d)$.
\item[B)] Denote by ${\zf}^{X_1\times X_2}$ the second fundamental form tensor
of $Y$ in $X_1\times X_2$. With ${\zf}^{X_1\times X_2}(A,B)=0={\zf}^{X_1\times X_2}(A,C)$
the Gauss equation yields
\begin{displaymath}
<R^Y_{A \; C}B,C>_Y \; = \; <R^Y_{A \; C}B,C>_{X_1\times X_2} \; = \;
<R^{X_1\times X_2}_{A \; C}B,C>_{X_1\times X_2}.
\end{displaymath}
Thus it follows from condition $v)$ that $B)$ is also bounded.
\item[C)] From $|K^{g_{X_{i}}}| \; < \; {\Lambda}_i$ and $iii)$ it follows from the Gauss equation that
\begin{displaymath}
\Big| {\zf}^i(v_i,v_i) {\zf}^i(w_i,w_i) \; - \; {\zf}^i(v_i,w_i) {\zf}^i(v_i,w_i) \Big|
\end{displaymath}
are bounded by some constants for all $g_{X_i}$-orthonormal systems $\{ v_i,w_i \}$,
which together with $ii)$ implies the existence of constants $\tilde{\epsilon}_i$ such
that 
\begin{displaymath}
\frac{{\zf}^i(v_i,v_i)}{g_{X_i}(v_i,v_i)} \; < \;  \tilde{\epsilon}_i .
\end{displaymath}
Hence we find as in $b)$ ${\zf}^Y(B,B) \; \le \; \tilde{\epsilon}$ for all $B \in TY$ 
tangent to $N$ of norm $||B||_{g_Y}=1$.
From 
\begin{displaymath}
\epsilon \; \le \; {\zf}^Y(B,B) \; \le \;  \tilde{\epsilon}
\end{displaymath}
for all $B \in TY$ tangent to $N$ of norm $||B||_{g_Y}=1$
one achieves 
\begin{displaymath}
{\epsilon}^2 \; \le \;
{\zf}^Y(B,B) {\zf}^Y(C,C) \; - \; {\zf}^Y(B,C) {\zf}^Y(B,C)
\; \le \; {\tilde{\epsilon}}^2
\end{displaymath}
for all $g_Y$-orthonormal systems $\{ B,C\}$. \\
Therefore it follows from the Gauss equation and condition $c)$ that the term $C)$ is also bounded.
\end{description}
\end{description}
Now a) - e) ensure that we can apply Lemma \ref{lemma}. Thus the validity of 
Proposition \ref{prophyprank} follows. 
\begin{flushright}
{\bf q.e.d.}
\end{flushright}

\begin{corollary} \label{corollary}
Let $(X_1,g_{X_1})$ and $(X_2,g_{X_2})$ be two Hadamard manifolds with pinched negative 
sectional curvature
$-b_i^2\le K^{X_i}\le -a_i^2 < 0$, $a_i,b_i \in \mathbb{R}$, $i=1,2$. Then their 
Riemannian product admits a Riemannian hypersurface that is 
bilipschitz to a Riemannian manifold $(Y,g^{\lambda}_Y)$ of negative sectional curvature 
bounded by $-b^2\le K^{g^{\lambda}_Y}\le -a^2 < 0$
for appropriate $a,b\in \mathbb{R}$.
\end{corollary} 
{\bf Proof of Corollary \ref{corollary}:} \\
Write the factors in horospherical coordinates and apply  Proposition \ref{prophyprank}. 
\begin{flushright}
{\bf q.e.d.}
\end{flushright}
Note that the additional boundary condition ensures that the horospheres' intrinsic curvature is also bounded (\cite{buka}). It is for that reason
that our method of prove requires this extra condition. \\


Since $(Y,g_Y^{\lambda)}$ is bilipschitz to $(Y,g_Y)$ it suffices to prove that 
the Riemannian embedding $(Y,g_Y)$ is bilipschitz to $(X,g_X)$. 

Denote the Riemannian length function on $Y$, that is induced by the metric $g_Y:=j^*g_X$ on $Y$, by $d_Y$ and that on
$X$, induced by the Riemannian product metric $g_X$, by $d_X$. With that we find the inequalities
\begin{displaymath}
d_X\Big( j(p),j(q)\Big) \; \stackrel{a)}{\le} \; d_Y\Big( p,q \Big)  \; \stackrel{b)}{\le} \; \Big( 2 \sqrt{2} + 2\Big) \cdot
d_X\Big( j(p),j(q)\Big) \hspace{0.5cm} \forall \; p,q \; \in \; Y.
\end{displaymath}  
While inequality a) merely is a consequence of the fact that $(Y,g_Y)$ is a Riemannian submanifold of $(X,g_X)$, 
inequality b) requires more attention. Note that the idea for the following construction is due to Brady and Farb (\cite{brfa}). \\
We denote the various canonical projections as follows:
\begin{displaymath} 
\begin{array}{lcl}
{\pi}_i: X \longrightarrow X_i  & \hspace{2cm} &     t : X_i \longrightarrow \mathbb{R}\\
\eta : Y \longrightarrow N & \hspace{2cm} & {\eta}_i : X \longrightarrow N_i 
\end{array} .
\end{displaymath}
Consider an arbitrary differentiable curve $c: [t_p,t_q] \longrightarrow X$ connecting $j(p)\in X$ with $j(q)\in X$. The idea
is to construct a curve $\tilde{c}:[\alpha ,\omega ]\longrightarrow Y$ that connects $p\in Y$ with $q\in Y$, whose Riemannian
length $L_Y(\tilde{c})$ in $(Y,g_Y)$ is bounded by a constant times the Riemannian length $L_X(c)$ of $c$
in $(X,g_X)$.\\
Therefore we consider the projections $c_i := {\pi}_i \circ c$ of $c$ to the factors $X_i$, $i=1,2$, that connect
$p_i:={\pi}_i(j(p))$ with $q_i:={\pi}_i(j(q))$. The further projections
$t_i \circ c_i$ are continuous, thus the set ${\cal K}:= \cup_{i=1,2}(t_i\circ c_i )([t_p,t_q]) \subset \mathbb{R}$ is compact
and therefore takes its maximum $(t_{i_0} \circ c_{i_0})(b)$ for some $b\in [t_p,t_q]$, $i_0 \in \{ 1,2 \} $.\\
The continuous and piecewise differentiable curve $\tilde{c}$ in $Y$ we are going to follow from $p$ to $q$ consists of three 
differentiable segments $v_1$, $\gamma$  and $v_2$ as follows:
\begin{itemize}
\item $v_1$ has constant projection $\eta \circ v_1$ to $N_1 \times N_2$, that is given through 
$\eta \circ v_1 \equiv ({\eta}_1(p_1), {\eta}_2(p_2))$, while its projection to $\mathbb{R}$ is
$t \circ v_1 = (t_{i_0} \circ c_{i_0})|_{[t_p,b]}$.
\item $\gamma$ is the curve keeping its projection to $\mathbb{R}$ constant: $(t\circ \gamma) \equiv (t_{i_0} \circ c_{i_0})(b)$,
while varying along $N_1 \times N_2$ with $\eta (\gamma ) = ({\eta}_1(c_1),{\eta}_2(c_2))$.
\item $v_2$ again has constant projection $\eta \circ v_2$ to $N_1 \times N_2$ that is 
$\eta \circ v_2 \equiv ({\eta}_1(q_1), {\eta}_2(q_2))$. Its projection to $\mathbb{R}$ is $t\circ v_2 = (t_{i_0} \circ c_{i_0})|_{[b,t_q]}$.
\end{itemize} 
The length of $\tilde{c}:= v_2 * \gamma * v_1$ is the sum
\begin{displaymath}
L_Y(\tilde{c}) \; = \; L_Y(v_1) \; + \;  L_Y(\gamma ) \; + \; L_Y(v_2).  
\end{displaymath}
From (\ref{metricg_Y}) it directly follows, that
\begin{eqnarray}
L_Y(v_m) & = & \sqrt{2} \cdot L_{(\mathbb{R},dt^2)} ((t_{i_0} \circ c_{i_0})|_{I_m}) \nonumber \\
& \le & \sqrt{2} \cdot L_{(\mathbb{R},dt^2)} (t_{i_0} \circ c_{i_0}) \nonumber \\   
& \le & \sqrt{2} \cdot L_X(c) , \label{absch1*} 
\end{eqnarray}
where $I_1=[t_p,b]$ and $I_2=[b,t_q]$. \\

Now we know that the natural diffeomorphisms ${\phi}^{tt'}_i : N^t_i \longrightarrow N^{t'}_i$ are length contracting for 
$t\le t'$. With that and the particular choice of $b$ it is 
\begin{eqnarray}
L_Y(\gamma ) & = & \int\limits_{t_p}^{t_q} \sqrt{g_Y({\gamma}',{\gamma}')} \; dt \nonumber \\
& = & \int\limits_{t_p}^{t_q} \Big[ g_Y{\Big|}_{N^{t_0}}\Big( (\eta \circ {\gamma})',(\eta \circ{\gamma})'\Big) 
{\Big]}^\frac{1}{2} \; dt \nonumber \\
& = & \int\limits_{t_p}^{t_q} \Big[ g_Y{\Big|}_{N_1^{t_0}}\Big( ({\eta}_1 \circ c)',({\eta}_1 \circ c)' \Big) \; + \; 
                                    g_Y{\Big|}_{N_2^{t_0}}\Big( ({\eta}_2 \circ c)',({\eta}_2 \circ c)'\Big) 
{\Big]}^{\frac{1}{2}} \; dt \nonumber \\
& \le & \int\limits_{t_p}^{t_q} \Big[ g_Y{\Big|}_{N_1^{t}}\Big( ({\eta}_1 \circ c)',({\eta}_1 \circ c)'\Big ) \; + \; 
                                    g_Y{\Big|}_{N_2^{t}}\Big( ({\eta}_2 \circ c)',({\eta}_2 \circ c)'\Big) 
{\Big]}^{\frac{1}{2}} \; dt \nonumber \\
&\le & 2 \cdot L_X(c). \label{absch2*}
\end{eqnarray}
where $t_0:=(t_{i_0}\circ c_{i_0})(b)$. \\
Thus with (\ref{absch1*}) and (\ref{absch2*}) we can conclude that for an arbitrary curve $c$ in $X$ connecting two points
$j(p)$, $j(q)$ $\in j(Y)\subset X$ there exists a curve $\tilde{c}$ in $Y$ connecting $p$ and $q$ with
\begin{displaymath}
L_Y(\tilde{c}) \; \le \; \Big( 2\sqrt{2} \; + \; 2 \Big) \; L_X(c).
\end{displaymath}
Thus the required inequality b) follows by the definitions of the Riemannian length functions. \\

Finally note that the generalization to products of finitely many
Hadamardmanifolds of pinched negative sectional curvature is straight
forward using once again that products of quasi-isometric maps are 
quasi-isometric.
\begin{flushright}
{\bf q.e.d.}
\end{flushright}


\section{Proof of Theorem \ref{theo-hyprank}}

Recall the definition of the hyperbolic rank as given in \cite{buysch}:
\begin{definition}
Let $X$ be a metric space. Then the hyperbolic rank, $rank_hX$, of $X$ is defined via
\begin{displaymath}
rank_h X \; := \; \sup\limits_{Y} \; dim \, {\partial}_{\infty} Y,
\end{displaymath}
where the supremum is taken over all locally compact $CAT(-1)$ Hadamard spaces $Y$ quasi-isometrically embedded 
into $X$ and $ {\partial}_{\infty} Y$ denotes the topological dimension of the boundary of $Y$.
\end{definition}
In \cite{buysch} the authors proved the subadditivity of this hyperbolic rank with respect to products
of Hadamard manifolds of pinched negative sectional curvature. This together with the following Lemma, 
that states the corresponding superadditivity, yields the proof of Theorem \ref{theo-hyprank}.
\begin{lemma} \label{hyprank}
The hyperbolic rank, $rank_h$, is superadditive with respect to Riemannian products
of Hadamard manifolds $(X_i,g_i)$ of pinched negative sectional curvature, i.e.,
for the Riemannian product $(X,g)$ of the $(X_i,g_i)$, $i=1,...,k$, one has
\begin{displaymath}
rank_h\Big( X,g_X\Big) \; \ge \; \sum\limits_{i=1}^k \, rank_h \Big( X_i,g_i\Big) .
\end{displaymath} 
\end{lemma}
{\bf Proof of Lemma \ref{hyprank}:} \\
From the Morse' quasi-isometric lemma it easily follows that $rank_h (X_i)=dim X_i -1$. 
From Theorem \ref{theo-hyprank} we further conclude $rank_h X \ge -k  + \sum\limits_{i=1}^k dim X_i $ and thus

\begin{displaymath}
rank_h X \; \ge \;  -k  + \sum\limits_{i=1}^k dim X_i \; = \; 
\sum\limits_{i=1}^k \, (dim X_i \, - \, 1)\; = \;
\sum\limits_{i=1}^k \, rank_h \Big( X_i,g_i\Big) .  
\end{displaymath}
\begin{flushright}
{\bf q.e.d.}
\end{flushright}   


{\footnotesize UNIVERSIT\"AT Z\"URICH, MATHEMATISCHES INSTITUT, WINTERTHURERSTRASSE 190, 
CH-8057 Z\"URICH, SWITZERLAND \\
E-mail addresses: $\;\;\;\;\;$ foertsch@math.unizh.ch $\;\;\;\;\;$ vschroed@math.unizh.ch}


\end{document}